\documentclass[11pt]{amsart}
\usepackage{latexsym}
\usepackage{amssymb}
\usepackage{amsmath}
\usepackage{graphicx}

\newtheorem{thm}{Theorem}[section]
\newtheorem{lem}[thm]{Lemma}
\newtheorem{cor}[thm]{Corollary}
\newtheorem{prop}[thm]{Proposition}

\theoremstyle{definition}
\newtheorem{defn}[thm]{Definition}
\newtheorem{example}[thm]{Example}
\newtheorem{rem}[thm]{Remark}
 
\numberwithin{equation}{section}

\newenvironment{Proof}[1]{\par\noindent{\bf
    Proof{#1}:}\quad}{}
\def\addcontentsline#1#2#3{\relax}

\long\outer\def\demo#1. #2\par{\medbreak\noindent {\bf#1.\enspace}
        {\rm#2}\par\ifdim\lastskip<\medskipamount\removelastskip
        \penalty55\medskip\fi}
\def\bdemo#1. #2\par{\medbreak\noindent {\bf#1.\enspace}{\rm#2}\par}
\def\edemo{\ifdim\lastskip<\medskipamount\removelastskip\penalty55\medskip\fi}
\def\eps{\varepsilon}

\def\g{{\bf g}}
\def\h{{\bf h}}
\def\x{{\bf x}}

\def\D{{\mathcal D}}

\def\N{{\bf N}}

\def\S{{\mathcal S}}
\def\Z{{\bf Z}}

\def\mod{\mathop{\rm mod}\nolimits}

\def\longto{\mathop{\longrightarrow}\limits}
\def\bfa{\mathop{a\kern-.5em{a}\kern-.5em{a}}\nolimits}
\def\notconverge{\mathop{\not\kern-.8em{\longto}}\limits}
\def\nequiv{\not\equiv}
\newcommand\XXX[1]{}
\def\bsq{\hfill\raisebox{-3mm}{$\blacksquare$}\;}

\newcommand\hsp[1]{\mbox{}\hspace{#1mm}} 
\newcommand\hspm[1]{\hsp{-#1}}

\begin{document}

\title[Nilpotent Groups are Round]
{Nilpotent Groups are Round}

\author{D. Berend}\address{Daniel Berend, Departments of Mathematics
and of Computer Science, Ben-Gurion University of the Negev,
Beer-Sheva 84105, Israel, \hspace{20mm} \mbox{}
\mbox{\hspace{3mm} Department of Mathematics,}
Rice University, Houston, TX 77251, USA}

\email{berend@math.bgu.ac.il}

\author{M. D. Boshernitzan}
\address{Michael D. Boshernitzan\\
Department of Mathematics \\
Rice University, Houston, TX 77251, USA}
\email{michael@math.rice.edu}

\begin{abstract} \XXX{Jan 24, 2006}  We define a notion of roundness for finite groups. 
Roughly speaking, a group is round if one can order its elements 
in a cycle in such a way that some natural summation operators map 
this cycle into new cycles containing all the elements 
of the group. Our main result is that this combinatorial property 
is equivalent to nilpotence. 
\end{abstract}

\maketitle

\setcounter{section}{0}
\section{Introduction and Main Results}
\setcounter{equation}{0}
\medskip

Given a finite group $G$ of order $n=|G|$, a {\it
cycle} in $G$ is a finite sequence of elements of $G$. A cycle
may be continued to a periodic sequence over $G$, and it will be
often more convenient to think of the cycle as this sequence.
A {\it 1-1 cycle} over $G$ is a cycle of length $n$, in which
every element of the group appears exactly once.

An $n$-length cycle $\g=(g_i)_{i=0}^{n-1}$ in $G$ is {\it $k$-round} if
for every $k$ integers $m_1,m_2,\ldots,m_k$ the cycle
$$\g_{m_1,m_2,\ldots,m_k}
    =\left(\prod_{j=1}^k g_{i+m_j}\right)_{i=0}^{n-1}=
        \Big(g_{i+m_1}g_{i+m_2}\ldots g_{i+m_k}\Big)_{i=0}^{n-1}$$
is a 1-1 cycle in $G$. (The addition of indices here and later 
is to be understood modulo the length $n$ of the cycle.)
Such a cycle is in particular a 1-1 cycle. (To see it, take
$m_1=m_2=\ldots=m_k=0$.)

The following definition plays a central role in the study initiated in this paper.

\begin{defn} \label{basic-definition}
\XXX{p1본본본본}      
\begin{enumerate}
\item{
$G$ is {\it $k$-round} if it admits a $k$-round cycle.
}
\item{
An $n$-length cycle in $G$ is {\it totally round} if it is $k$-round for
every positive integer $k$ with $(k,n)=1$.
}
\item{
$G$ is {\it round} if it admits a totally round cycle.
}
\end{enumerate}
\end{defn}

Note that $G$ cannot possibly
be $k$-round if $(k,n)>1$. In fact, if $p$ is any
prime divisor of $(k,n)$ and $\g=(g_i)_{i=0}^{n-1}$ is any 1-1
cycle, then the cycle $\g^k=(g_i^k)_{i=0}^{n-1}$ contains
the identity element $1\in G$ at least $p$ times. This explains
the constraint $(k,n)=1$ in Definition \ref{basic-definition}
and Theorem \ref{characterization-of-round-groups} below.

Also, the property of a group being round appears 
to be stronger than the property of being $k$-round
for every $k$ with $(k,n)=1$. Nevertheless, the implication
$(2)\hspm2\implies\hspm2(3)$ in 
Theorem~\ref{characterization-of-round-groups}
will show that the two are in fact equivalent. \XXX{3본본본본}

Our main result is the equivalence of the last two conditions in the
following theorem.

\begin{thm} \label{characterization-of-round-groups}
Let $G$ be a finite group of order $n$. The following conditions
are equivalent:
\begin{enumerate}
\item
$G$ is $k$-round for some $k>n^2$.
\item
$G$ is $k$-round for every $k$ with $(k,n)=1$.
\item
$G$ is round.
\item
$G$ is nilpotent.
\end{enumerate}
\end{thm}

For additional equivalent conditions (for a finite group to be round),
see Theorem \ref{strongly-round}
and Remark \ref{rem:more_round} below.
Recall that one of the equivalent conditions for a finite group to be
nilpotent is to be isomorphic to a direct product of $p$-groups. 
(For this and other basic results used throughout the paper, 
we refer to any standard text, say \cite{herstein}.) \XXX{4본본본본} 

Theorem \ref{characterization-of-round-groups} implies that the
family of round groups satisfies some closure properties.

\begin{cor}
Any subgroup and any quotient group of a round group is round as well.
\end{cor}

In a similar vein, if $G$ is $k$-round and $l|k$ ($l$ divides $k$), then $G$ is $l$-round
as well. Moreover, if $\g$ is a $k$-round cycle, then it is also
$l$-round. In fact, let $\g'$ be any product of $l$ rotates of $\g$. Then
$\g'$ is obviously $k/l$-round, and, in particular, it is 1-1.  \XXX{5본본본본} 

\begin{example} \label{cyclic}
Any finite cyclic group $\Z_n$ is immediately seen to be round, as the
``arithmetic progression" cycle
$\g=\left( i\right)_{i=0}^{n-1}$ is totally round. Indeed, every product $\g'$
of $k$ rotates of $\g$ is itself an arithmetic progression whose
difference is $k$; thus $\g'$ is a 1-1 cycle if $(k,n)=1$.
The fact that other finite nilpotent groups, even abelian groups such
as $\Z_3\times \Z_3$, are round is less obvious. \XXX{6본본본} 
\end{example}

\begin{rem} \label{frobenius-number}
The first condition in Theorem \ref{characterization-of-round-groups}, 
which appears weaker than the next two, 
may be replaced by the even weaker condition
\begin{equation}\label{weakest}
G \text{ is } k\text{-round for some } k>\Phi'(n),\hspace{55mm}, \mbox{}
\end{equation}
\noindent where $\Phi'(n)$ is the maximal integer which
does not belong to the additive semigroup generated by
the (relatively prime) numbers
\begin{equation}\label{bi}
b_i=n/p_i^{e^i},\ 1\le i\le r,
\end{equation}
$n=p_1^{e_1} p_2^{e_2}\ldots p_r^{e_r}$
being the prime power factorization of $n$. 
(This is the essence of the proof of implication
$(1)\!\!\implies\!\!(4)$ in
Theorem~\ref{characterization-of-round-groups}.) 
Condition~\eqref{weakest} indeed implies condition (1)
in Theorem \ref{characterization-of-round-groups}; see \eqref{ineq}.
\XXX{7본본본볗

More generally, for a nonempty subset $S\subseteq\N=\{1,2,3,\ldots\}$ \ such that
$\gcd(S)=1$, denote by $\Sigma(S)$ the additive semigroup
generated by $S$, and put $\Sigma'(S)=\N-\Sigma(S).$
Thus $S\subseteq\Sigma(S)\subseteq\N$, and
$\Sigma'(S)$ is the set of positive integers which cannot 
be represented as finite sums of elements of $S$.
The number
$$\Phi(S)=\sup (\Sigma'(S)\cup \{ -1\})$$
is the so-called 
 {\em Frobenius number} of $S$. It is easily seen that  
 \XXX{8, 9본본본본본본볗
       $$
           \Phi(S)=-1 \iff 1\in S; \qquad \qquad \Phi(S)<\infty \iff \gcd (S)=1.
        $$
  In our case, denoting $B=\{b_i\mid 1\leq i\leq r\}$ (see \eqref{bi}),
we observe that $\gcd(B)=1$ and $\Phi'(n)=\Phi(B)$.

For $|S|=2$, there exists a simple formula for the 
Frobenius number $\Phi(S)$, which
goes back at least as far as Sylvester \cite{sylvester}: If
$(x,y)=1$, then
$
\Phi(\{x, y\})=xy-x-y.
$
However, there exists no such formula for $|S|\geq 3$, 
nor should such a formula be expected. (See, for example, 
\cite{bak}, \cite{beck-zacks} and the references there for more information 
regarding the so-called linear diophantine problem of Frobenius.) 
Nevertheless, there are various upper bounds: 
\cite{erdos-graham},\cite{vitek},\cite{selmer}.
Employing these bounds, we easily see that 
\begin{equation} \label{ineq}
\Phi'(n)=\Phi({\bf B})\leq n^2.
\end{equation}
(Actually, better estimates are possible for $\Phi'(n)$, but we care here about 
the phenomenon rather than the exact bound.) 
\XXX{10본본본본본본볗
\end{rem}

\begin{rem}\label{rem:more_round}
It is possible to suggest a stronger notion of roundness.
Let a cycle $\g=(g_i)_{i=0}^{n-1}$ over $G$ be
{\it strongly $k$-round} if for every number $l$ of integers
$m_1,m_2,\ldots,m_l$ and $\eps_1,\eps_2,\ldots,\eps_l\in\{-1,1\}$
with $\sum_{j=1}^l \eps_j=k$, the cycle
$\left(\prod_{j=1}^l g_{i+m_j}^{\eps_j}\right)_{i=0}^{n-1}$
is 1-1. 
Observe that, if $k_1\equiv k_2 (\mod n)$, then a cycle is strongly 
$k_1$-round if and only if it is strongly $k_2$-round.
Indeed, multiplying a product of rotates
of a cycle $\g$ by $\g^n$ (or $\g^{-n}$) does not alter
the product.
We conclude that the existence of a strongly $k$-round cycle 
for some $k$ implies the existence of such a cycle for some 
$k>n^2$. This cycle is in particular $k$-round, and, by Theorem
\ref{characterization-of-round-groups}, the group $G$ is round.
On the other hand,
a totally round cycle is clearly strongly $k$-round for every~$k$.
Consequently, a group admits a strongly $k$-round cycle for some $k$
if and only if it is round.
\end{rem}

The above discussion is summarized in the 
following theorem.
\XXX{11본본본본}
\pagebreak[4]

\begin{thm} \label{strongly-round}
Let $G$ be a finite group of order $n$. The following four
conditions are equivalent:
\begin{enumerate}
\item
$G$ admits a strongly $1$-round cycle.
\item
$G$ admits a strongly $k$-round cycle for some $k$.
\item
$G$ admits a strongly $k$-round cycle for every $k$ with $(k,n)=1$.
\item
$G$ is round.
\end{enumerate}
\end{thm}

While the proof of Theorem \ref{characterization-of-round-groups} gives
a way of constructing totally round cycles for nilpotent groups, it
does not provide a way of deciding whether a given cycle is such. The
following proposition gives an algorithm to that effect.  
\XXX{12본본본본볗

\begin{prop} \label{algorithm-for total-roundness}
There exists an effective constant $K=K(n)$ such that, if a cycle $\g$ of
length $n$ over $G$ is $k$-round for every $k\le K$ with $(k,n)=1$,
then it is totally round.
\end{prop}

The following proposition demonstrates that a non-nilpotent group may
still be $k$-round for some $k$'s. The simplest example with $k=2$ is
given by the family of non-abelian groups whose order is a product of
two distinct odd primes. Recall that, for odd primes $p,q$
with $p<q$, a group of order $pq$ exists (and is unique) if and only
if $q\equiv 1 (\mod p)$.
\XXX{13본본본본볗

\begin{prop} \label{pq-group}
 Let $G$ be a non-abelian group of order $pq$, where $p,q$ are
odd primes with $p<q$. Then $G$ is $2$-round.
\end{prop}

 Nevertheless, some groups
are {\it unround} -- do not admit a $k$-round
cycle for any $k>1$. The following theorem presents two such families
of groups.

\begin{thm} \label{unround-groups}
The following groups are unround:
\begin{enumerate}
\item Dihedral groups of order divisible by 3.
\item The symmetric group $S_l$ for every $l\ge 3$. More generally,
all sufficiently large almost simple groups with the exception of
the groups $^2 B_2(q)$ (cf. \cite{hurley-rudvalis}).
\end{enumerate}
\end{thm}

One of the motivations for this study (which is a special case of
Example \ref{cyclic}) was
the observation that a cyclic group of odd order admits
a 1-1 cycle, whose
product with its rotate (by a single position) is again such. 

More precisely, define an
operator $\S$ (for sum) on the set of all cycles
$\g=(g_i)_{i=0}^{n-1}$ over $G$ by:
\XXX{14} 
$$\S(\g)=\g_{0,1}=(g_i g_{i+1})_{i=0}^{n-1}, \qquad \g\in G^n.$$
The following proposition follows in a straightforward manner from
Theorem \ref{characterization-of-round-groups} and shows that in many
cases we may apply $\S$
over and over again, obtaining each time a 1-1 cycle.

\begin{prop}
A nilpotent group of odd order admits a cycle $\g$ such that
$\S^m(\g)$ is a 1-1 cycle for every positive integer $m$.
\XXX{15본본본본본본본} 
\end{prop}

In view of the above proposition, it is natural to ask
 whether there exists an analogue if we replace addition by subtraction. 
 Of course, starting with a
1-1 cycle $\g=(g_i)_{i=0}^{n-1}$ in $G$, the cycle consisting of
``differences" of consecutive entries of $\g$, namely $(g_i
g_{i-1}^{-1})_{i=0}^{n-1}$, does not contain the unit element of $G$, and
thus cannot possibly be 1-1. 
Thus, we shall consider now cycles of
arbitrary lengths $L$, where $n|L$. Such a cycle is {\it balanced} if it
contains the same number of occurrences of each element of $G$. 
Define a transformation $\D$ 
on the set of all cycles $\g=(g_i)_{i=0}^{L-1}$
of length $L$ over $G$ by:
\XXX{16본본본본본본볗 
$$\D(\g)=(g_i g_{i+1}^{-1})_{i=0}^{L-1}.$$
For an integer $l\geq 1$,
a cycle $\g$ is $\D^l${\em -balanced} if $\D^r(\g)$ is
balanced for every integer $0\leq r\leq l$. 
A cycle $\g$ is $\D^\infty${\em -balanced} if $\D^r(\g)$ is
balanced for every $r\geq 1$. Finally,
for $1\leq s\leq\infty$, a group $G$ is $\D^s${\em -balanced}
if it admits a $\D^s${-balanced} cycle.

The following result will be
shown to follow relatively easily from Theorem~\ref
{characterization-of-round-groups}.

\begin{thm} \label{balanced-cycles}
Every finite nilpotent group of odd order is $\D^\infty$-balanced.
\end{thm}

This result should be contrasted with

\begin{thm} \label{balanced-cycles2}
Every finite group is $\D^r$-balanced for all integers $r\geq 1$.
\end{thm}

\begin{rem}
For non-abelian groups one may consider four different versions
for the transformation $\D$. Namely, starting from
$\g=(g_i)_{i=0}^{L-1}$, one may let the $i$-th entry of
$\D(\g)$ be either $g_i g_{i+1}^{-1}$, as done above, or
$g_{i+1}^{-1} g_i$, or $g_i^{-1} g_{i+1}$, or $g_{i+1} g_i^{-1}$.
It turns out that both Theorems \ref{balanced-cycles} 
and \ref{balanced-cycles2} hold for each of
these variations. Moreover, the cycle $\g$ we construct actually 
has the property that we may apply at each stage any of the four
variations (of the difference operator)
and still get only balanced cycles.
\XXX{17}
\end{rem}

We conjecture that the group $\Z_2=\{0,1\}$ is not $\D^\infty$-balanced.
The conjecture would imply that no solvable finite group of even order is
 \mbox{$\D^\infty$-balanced}.

The authors wish to express their gratitude to V. Lev for referring them
to his paper \cite{seva}, which served as
the initial motivation for considering
the questions studied in this paper.
\vspace{0mm}

\section{Proofs}
\medskip
Before proving Theorem \ref{characterization-of-round-groups}, we 
present several auxiliary results.

\begin{lem} \label{normal-subgroup}
If $G$ has a subgroup $H$ contained in the center $Z(G)\subseteq G$,
such that both $H$ and $G/H$
are $k$-round, then $G$ is $k$-round as well.
\end{lem}

\begin{lem} \label{sum-of-groups}
If $G_1$ and $G_2$ are both $k$-round, then $G_1\times G_2$
is $k$-round as well.
\end{lem}

Obviously, the following lemma contains both
Lemmas \ref{normal-subgroup} and \ref{sum-of-groups} as
special cases, so that we shall prove only it instead.
\XXX{18본본본본본} 

\begin{lem}\label{general-extension}
If $G$ has a normal subgroup $H$, each coset of which contains
an element commuting with all elements of $H$, and if both $H$
and $G/H$ are $k$-round, then $G$ is $k$-round as well.
\end{lem}

\begin{rem} \label{universality}
In fact, we prove more than is stated. Namely, the construction of a
$k$-round cycle in $G$, based on such cycles in $H$ and $G/H$, does not
depend on $k$. Hence, if we start with cycles in these two groups which
are $k$-round for several values of $k$, then so are the resulting
cycles in $G$. In particular, if both $H$ and $G/H$ are round, then so
is $G$.
\XXX{19}
\end{rem}

\begin{Proof}{\ of Lemma \ref{general-extension}}
 Let  $H=\{h_0,h_1,\ldots,h_{s-1}\}$  and  let
$x_0, x_1,\ldots,x_{t-1}$
be representatives of all cosets of $H$, such that each $x_j$ commutes
with $H$. Rearranging the $h_i$'s and $x_j$'s if necessary, we may
assume that the cycle   $\h=(h_i)_{i=0}^{s-1}$ is  $k$-round
in $H$ and that the cycle $\x=(Hx_j)_{j=0}^{t-1}$ is
$k$-round in~$G/H$. Consider the cycle
$$\g=(h_0 x_0,\ldots,h_0 x_{t-1},h_1 x_0,\ldots,h_1 x_{t-1},\ldots,
      h_{s-1} x_0,\ldots,h_{s-1} x_{t-1})$$
in $G$. (Formally, if $0\le i\le n-1$, write $i=at+b$ with 
\mbox{$0\le a\le s-1$},  \mbox{$0\le b\le t-1$},  and then put 
$g_i=h_a x_b=h_{\lfloor i/t\rfloor} x_{i\,\mod t}$.)
We claim that $\g$ is $k$-round. Indeed, suppose it
is not. Let $m_1,m_2,\ldots,m_k$ be integers in
$\{0,1,\ldots,n-1\}$ such that 
\begin{equation} \label{basic-equation}
\prod_{j=1}^k g_{i+m_j}=\prod_{j=1}^k g_{i'+m_j}
\end{equation}
for some $i,i'$
with $0\le i,i'\le n-1$. Projecting to $G/H$, we obtain
$$\prod_{j=1}^k Hg_{i+m_j}=\prod_{j=1}^k Hg_{i'+m_j},$$
that is
\begin{equation} \label{projected-equation}
H\prod_{j=1}^k x_{i+m_j\,(\mod t)}=H\prod_{j=1}^k x_{i'+m_j\,(\mod t)}.
\end{equation}

Since $\x$ is a $k$-round cycle in $G/H$, (\ref{projected-equation})
implies $i\equiv i' (\mod t)$. Hence (\ref{basic-equation}) may be
rewritten in the form
\XXX{20} 
$$\prod_{j=1}^k h_{\lfloor (i+m_j)/t\rfloor\,\mod s}\; x_{i+m_j\,(\mod t)}
  =\prod_{j=1}^k h_{\lfloor (i'+m_j)/t\rfloor\,\mod s}\; x_{i'+m_j\,(\mod t)}.$$
As the $x_i$'s commute with $H$, this yields
$$\prod_{j=1}^k h_{\lfloor (i+m_j)/t\rfloor\,\mod s}
   =\prod_{j=1}^k h_{\lfloor (i'+m_j)/t\rfloor\,\mod s}.$$
Since $\h$ is a $k$-round cycle in $H$, this implies that
$\lfloor i/t\rfloor\,\mod s=\lfloor i'/t\rfloor\,\mod s$,
and therefore $i=i'$. This proves the proposition. \bsq
\end{Proof}
\medskip

 Our main tool for proving that a group is not $k$-round for a
certain~$k$ is provided by the following proposition. For a finite group $G$
and a positive integer $l$, denote by $R_l(G)$ the number of solutions of the equation $x^l=1$
in~$G$.
\XXX{21}

\begin{prop} \label{condition-for-non-k-roundness}
If
$$R_{l_1}(G)R_{l_2}(G)\ldots R_{l_s}(G)>n^{s-1}$$
for some positive integers
$s,l_1,l_2,\ldots,l_s$, then $G$ is not $k$-round for any $k$ of the form
$k=c_1 l_1 + c_2 l_2 + \ldots + c_s l_s$, where $c_1,c_2,\ldots,c_s$
are non-negative integers.
\end{prop}

 Obviously, one would usually apply the proposition with all $l_i$'s
being divisors of $n$ and strictly positive $c_i$'s. We may view the
proposition as a strengthening of our former observation (following
Definition \ref{basic-definition}) that $G$
cannot be $k$-round if $(k,n)>1$.
 \XXX{22}
\begin{Proof}{}
Let $\g=(g_i)_{i=0}^{n-1}$ be any cycle, and let $k$ be as in the
proposition. We have to show that $\g$ is not $k$-round.

Consider the cycles $\g^{(j)}=\g^{c_j l_j},\ 1\le j\le s$. According to our
assumptions, the multiplicity $M_j$ of the identity element $1\in G$
in $\g^{(j)}$ is at least $R_{l_j}(G)$. Now consider all cycles of the form
$\g_{m_1,\ldots,m_s}=\left(\prod_{j=1}^s \g^{(j)}_{i+m_j}\right)_{i=0}^{n-1}$ as
$(m_1,\ldots,m_s)$ runs through all $s$-tuples of integers between $0$
and $n-1$. We have $n^s$ cycles, containing between them altogether at
least $n\prod_{j=1}^s M_j$ occurrences of 1. Hence at least one of
these cycles contains at least
$$n(\prod_{j=1}^s M_j)/n^s\ge (\prod_{j=1}^s R_{l_j}(G))/n^{s-1}>1$$
occurrences of 1. Hence $\g$ is not $k$-round.
\XXX{23}  \bsq
\end{Proof}

\medskip
 In the course of the proof we shall also use the following result of
Frobenius~\cite{frobenius}.

\begin{thm} \label{frobenius}
Let $G$ be a finite group of order $n$ and let $l|n$. Then $l|R_l(G)$.
\end{thm}

The case where $R_l(G)=l$ is of particular interest. The following
theorem, conjectured by Frobenius \cite{frobenius}, waited almost
a century for a proof.

\begin{thm} \cite{iiyori}\label{iiyori-frobenius}
In the setup of Theorem \ref{frobenius}, if $R_l(G)=l$, then
the set $\{x\in G:x^l=1\}$ is a normal subgroup of $G$.
\end{thm}

\begin{Proof}{\ of Theorem \ref{characterization-of-round-groups}} 

$(2) \Rightarrow (1)$:
Trivial.

$(3) \Rightarrow (2)$:
Trivial.  
\XXX{24본본본본본} 

$(4) \Rightarrow (3)$:
Since a $p$-group has a non-trivial
center, Example \ref{cyclic}, Lemma \ref{normal-subgroup} and
Remark \ref{universality} imply that every $p$-group is round. 
Lemma \ref{sum-of-groups} (again with Remark \ref{universality})
then implies that every nilpotent group is round.

$(1) \Rightarrow (4)$:
Suppose $G$ is $k$-round for some $k>n^2$. As explained in Remark
\ref{frobenius-number}, this implies that we may represent
$k$ in the form
$$k=n\left(\frac{c_1}{p_1^{e_1}}+\frac{c_2}{p_2^{e_2}}+\ldots
   +\frac{c_r}{p_r^{e_r}}\right)
    $$
for suitable positive integers $c_1,c_2,\ldots,c_r$, where
$$n=p_1^{e_1} p_2^{e_2}\ldots p_r^{e_r}$$
is the prime power factorization of $n$.
We have to show that $G$ is nilpotent.
Consider the numbers $M_j=R_{n/p_j^{e_j}}(G)$ of solutions of the
equations \mbox{$x^{n/p_j^{e_j}}=1$ in $G$.}
According to Theorem \ref{frobenius}, applied with $l=n/p_j^{e_j}$,
each $M_j$ is a multiple of $n/p_j^{e_j}$, and in particular 
$M_j\ge n/p_j^{e_j}$. If for some $j$ this inequality is strict, 
then $\prod_{j=1}^r M_j>n^{r-1}$,
so that, by Proposition \ref{condition-for-non-k-roundness}, $G$ is not
$k$-round. Thus $M_j=n/p_j^{e_j}$ for each $j$.
According to Theorem \ref{iiyori-frobenius},
this implies that each of the sets
$$H_j=\{g\in G:g^{n/p_j^{e_j}}=1\}$$
is a subgroup of $G$. Consider the subgroups:
$$H'_j=\bigcap_{l\ne j} H_l, \qquad j=1,2,\ldots,r.$$

We observe that $|H'_j|\big|(n/p_l^{e_l})$ for each $l\ne j$, 
and therefore $|H'_j|\big|p_j^{e_j}$. 
On the other hand, $h\in H'_j$ if and only if
$h^{p_j^{e_j}}=1$, so that $H'_j$ is the union of all
$p_j$-Sylow subgroups of $G$. 
\XXX{26본본본본본볗

It follows that $G$ contains a unique $p_j$-Sylow subgroup for each $j$. 
 Hence $G$ is the direct
product of its Sylow subgroups, which means that it is nilpotent. 
This completes the proof. \bsq
\end{Proof}

\medskip

\begin{Proof}{\ of Proposition \ref{algorithm-for total-roundness}}
Take $K=(n+1)!$. Given an 
\mbox{$n$-length} cycle $\g$, which is $k$-round
for every $k\le K$ with $(k,n)=1$, we have to show that it is
totally round. In fact, suppose it is not totally round. Let $k_0$
be the minimal~$k$ with $(k,n)=1$ for which $\g$ fails to be
$k$-round. Take $m_1,m_2,\ldots,m_k$ for which
$\g_{m_1,m_2,\ldots,m_k}$ is not 1-1. Put $k'=k_0 \mod n$. Consider
the cycles
$$\g_{m_1,m_2,\ldots,m_{k'}},\ \g_{m_1,m_2,\ldots,m_{k'+n}},\
\g_{m_1,m_2,\ldots,m_{k'+2n}},\ \ldots,\
\g_{m_1,m_2,\ldots,m_{k'+n!\cdot n}}.$$
According to our assumptions, all these are 1-1. Hence two of them
are identical, say
$$\g_{m_1,m_2,\ldots,m_{k'+jn}}=\g_{m_1,m_2,\ldots,m_{k'+j'n}}.$$
It follows that
$$\g_{m_1,m_2,\ldots,m_k}=\g_{m_1,m_2,\ldots,m_{k'+jn},m_{k'+j'n+1},
   \ldots,m_{k_0}}.$$
The cycle on the right-hand side is a product of less than $k_0$
rotates of $\g$, and thus should be 1-1, which yields a
contradiction. \bsq
\end{Proof}
\medskip

\begin{Proof}{\ of Proposition \ref{pq-group}}

We have to construct a 2-round cycle in $G$. Denote
$r=\frac{q-1}{p}$, and let $u$ be an integer such that
\XXX{27본본본본본볗 
\begin{equation} \label{properties-of-u}
u\nequiv 1 (\mod q),\qquad u^p\equiv 1 (\mod q).
\end{equation}
The group $G$ is
generated by $\{x,y\}$, satisfying the relations
$$\begin{array}{ll}
   & x^q=y^p=1,\\
   & yx = x^u y.
  \end{array}$$
The subgroup $H=\{1,x,x^2,\ldots,x^{q-1}\}$ is normal in $G$. 
Consider the cycle
$$\g=(1,y,\ldots,y^{p-1},x,xy,\ldots,xy^{p-1},\ldots,
      x^{q-1},x^{q-1}y,\ldots,x^{q-1}y^{p-1}).$$
We claim that $\g$ is $2$-round. To this end, we have to show
that, given an integer $m$, the cycle
$\g_{0,m}=\left(g_i g_{i+m}\right)_{i=0}^{pq-1}$ is a 1-1 cycle.
Indeed, suppose
\begin{equation} \label{supposed-equality}
g_i g_{i+m}=g_i' g_{i'+m}
\end{equation}
for some $i,i'$ with $0\le i,i'\le pq-1$. Similarly to the proof of
Proposition \ref{general-extension}, this can be shown to imply,
since $G/H$ is isomorphic to $\Z_p$, that  \mbox{$i\equiv
i' (\mod p)$.} Thus for some $a_1,a_2,b_1,b_2,c$ we may rewrite
(\ref{supposed-equality}) in the form
$$x^{a_1} y^{b_1} x^{a_2} y^{b_2} = x^{a_1+c} y^{b_1} x^{a_2+c}
y^{b_2},$$
which yields
\XXX{28본본본본본볗 
$$x^{a_1+a_2 u^{b_2}}=x^{a_1+c+(a_2+c) u^{b_2}}\,.$$
Consequently
$$c\left(1+u^{b_2}\right)\equiv 0 (\mod q),$$
which implies, since by (\ref{properties-of-u}) the second factor on the
left-hand side cannot vanish modulo~$q$, that $c\equiv 0 (\mod q)$.
Hence in (\ref{supposed-equality}) we have $i\equiv i' (\mod pq)$, so
that $\g$ is indeed $2$-round. \bsq
\end{Proof}
\medskip

\begin{Proof}{\ of Theorem \ref{unround-groups}}
We use Proposition \ref{condition-for-non-k-roundness}.

For a dihedral
group $D_l$ with $3|l$, we need to show that the group is not $k$-round
for odd $k>3$. Such a $k$ may be written in the form $c\cdot 2+3$, and
consequently the inequality
\XXX{29본본본본본볗 
$$R_2(D_l) R_3(D_l) = (l+1)\cdot 3>2l=|G|$$
proves that the group is not $k$-round.

For an almost simple group $G$, other than $^2B_2(q)$,
we have in view of \cite[Sec. 4]
{liebeck-shalev-1} and \cite[Prop. 3.1,3.2]{liebeck-shalev-2}
$$R_2(G) R_3(G) \ge c |G|^{1/2}\cdot c |G|^{3/5}>|G|$$
for an appropriate constant $c$. As for the dihedral group, this implies
that $G$ is not $k$-round for any $k$.

In the specific case of $S_l$, there is a lot of information regarding
the numbers $R_j(G)$ (cf. \cite{herrera} and the references there). In
particular, denoting $R_{2,l}=R_2(S_l)$, it is easy to prove the
recurrence
$$R_{2,l}=R_{2,l-1}+(l-1)R_{2,l-2},$$
from which it follows \cite{chowla-herstein-moore} that
$R_{2,l}/R_{2,l-1}>\sqrt{l}$, so that $R_{2,l}>\sqrt{l!}$ for $l\ge 2$.
Similarly, denoting $R_{3,l}=R_3(S_l)$, it is easy to prove the recurrence
$$R_{3,l}=R_{3,l-1}+(l-1)(l-2)R_{3,l-3},$$
which implies by an easy induction that $R_{3,l}>l!^{2/3}$ for $l\ge
3$. In particular
$$R_{2,l} R_{3,l} >l!, \qquad l\ge 3,$$
so that $S_l$, for $l\ge 3$, is not $k$-round for any $k$.
\XXX{30본본본본본볗
\bsq

\end{Proof}
\medskip

\begin{Proof}{\ of Theorem \ref{balanced-cycles}}
Let $(g_i)_{i=0}^{n-1}$ be a totally round cycle over $G$. We shall show
that the $2n$-length cycle
$$\g'=(g_0,g_0^{-1},g_1,g_1^{-1},g_2,g_2{-1},\ldots,g_{n-1},g_{n-1}^{-1})$$
is a $\D^\infty$-balanced cycle. In fact, it is easy to check
that, for each $r$, the cycle $\D^r(\g')$ consists of a merge
of two 1-1 cycles, as follows. The length-$n$ cycle consisting of
the entries at the places $0,2,4,\ldots,2(n-1)$ is obtained by
multiplying $2^r$ rotates of $\g$. The length-$n$ cycle consisting
of all other entries is obtained by inverting all entries in a
product of $2^r$ rotates of $\g$. Since $\g$ is in particular
$2^r$-round, each of these subcycles is a 1-1 cycle, and therefore
the whole cycle is balanced. \bsq

\end{Proof}
\medskip

\begin{Proof}{\ of Theorem \ref{balanced-cycles2}}
Let $|G|=n$. There exists a 1-1 cycle $\g=(g_i)_{i=0}^{L-1}$,
with $L=n^{r+1}$, such that each of the $L$ possible $(r+1)$-blocks of elements
in $G$ appears in $\g$ exactly once. (Such cycles exist; these are 
the so-called
{\em complete cycles of order} $r+1$ in $G$, or {\em De Bruijn sequences} --
see \cite{dB1}, \cite[pp.91--99]{Hall}). 
It is straightforwardly verified that such
cycles $\g$ must be $\D^r$-balanced. \bsq
\end{Proof}


\bigskip

\end{document}